\pdfoutput=1
\IfFileExists{./prepreamble-amspreprint.sty}{\RequirePackage[packages,theorems,changes]{./prepreamble-amspreprint}}{}

\makeatletter
\IfFileExists{./scoop-latex/scoop-packages.sty}{\providecommand*{\input@path}{}\edef\input@path{{./scoop-latex/}\input@path}}{}
\@namedef{ver@minted.sty}{}
\@namedef{opt@minted.sty}{}
\makeatother
\PassOptionsToPackage{style=authoryear-comp,backend=biber}{biblatex}    

\documentclass[english]{amsart}

\RequirePackage{biblatex}
\RequirePackage{ltxcmds}
\IfFileExists{./preamble-amspreprint.sty}{\RequirePackage[packages,theorems,changes]{./preamble-amspreprint}}{}

\usepackage{scoop-local}

\IfFileExists{./postpreamble-amspreprint.sty}{\RequirePackage[packages,theorems,changes]{./postpreamble-amspreprint}}{}

\makeatletter
\@ifpackageloaded{changes}{
\definechangesauthor[name = {Alexandra Bünger}, color = {red!80!black}]{AB}
\definechangesauthor[name = {Roland Herzog}, color = {blue!80!black}]{RH}
\definechangesauthor[name = {Andreas Naumann}, color = {orange!80!black}]{AN}
\definechangesauthor[name = {Martin Stoll}, color = {green!70!black}]{MS}
}{}
\makeatother        

\makeatletter
\@ifpackageloaded{hyperref}{%
	\hypersetup{
		pdftitle = {Uncertainty Propagation of Initial Conditions in Thermal Models},
		pdfauthor = {Alexandra Bünger, Roland Herzog, Andreas Naumann, Martin Stoll},
		pdfkeywords = {parameter uncertainty, heat equation, PDE systems, sensor placement},
		pdfcreator = {Created using the Scoop Template Engine version 1.1.1.}
	}
}{
	\pdfinfo{
		/Title (Uncertainty Propagation of Initial Conditions in Thermal Models)
		/Author (Alexandra Bünger, Roland Herzog, Andreas Naumann, Martin Stoll)
		/Subject ()
		/Keywords (parameter uncertainty, heat equation, PDE systems, sensor placement)
		/Creator (Created using the Scoop Template Engine version 1.1.1.)
	}
}
\makeatother

\title[Uncertainty Propagation of Initial Conditions in Thermal Models]{Uncertainty Propagation of Initial Conditions in Thermal Models}

\author{Alexandra Bünger}
\address[A. Bünger]{Technische Universität Chemnitz, Faculty of Mathematics, Professorship Scientific Computing, 09107 Chemnitz, Germany}
\email{alexandra.buenger@mathematik.tu-chemnitz.de}
\urladdr{https://www.tu-chemnitz.de/mathematik/wire/people/buenger.php}

\author{Roland Herzog\orcidlink{0000-0003-2164-6575}}
\address[R. Herzog]{Interdisciplinary Center for Scientific Computing, Heidelberg University, 69120 Heidelberg, Germany}
\email{roland.herzog@iwr.uni-heidelberg.de}
\urladdr{https://scoop.iwr.uni-heidelberg.de/team/roland-herzog}

\author{Andreas Naumann\orcidlink{0000-0002-9287-0315}}
\address[A. Naumann]{Interdisciplinary Center for Scientific Computing, Heidelberg University, 69120 Heidelberg, Germany}
\email{andreas.naumann@iwr.uni-heidelberg.de}
\urladdr{https://scoop.iwr.uni-heidelberg.de/team/andreas-naumann}

\author{Martin Stoll\orcidlink{0000-0003-0951-4756}}
\address[M. Stoll]{Technische Universität Chemnitz, Faculty of Mathematics, Professorship Scientific Computing, 09107 Chemnitz, Germany}
\email{martin.stoll@mathematik.tu-chemnitz.de}
\urladdr{https://www.tu-chemnitz.de/mathematik/wire/prof.php}

\thanks{This work was supported by a DFG grant within the Special Research Program SFB/Transregio~96 (\emph{Thermo-energetische Gestaltung von Werkzeugmaschinen}) (DFG project \href{https://gepris.dfg.de/gepris/projekt/174223256}{1742243256}), which is gratefully acknowledged.}

\date{}

\dedicatory{}

\begin{document}

% Insert the abstract.
\begin{abstract}
The operation of machine tools often demands a highly accurate knowledge of the tool center point's (TCP) position. 
The displacement of the TCP over time can be inferred from thermal models, which comprise a set of geometrically coupled heat equations.
Each of these equations represents the temperature in part of the machine, and they are often formulated on complicated geometries.

The accuracy of the TCP prediction depends highly on the accuracy of the model parameters, such as heat exchange parameters, and the initial temperature.
Thus it is of utmost interest to determine the influence of these parameters on the TCP displacement prediction.
In turn, the accuracy of the parameter estimate is essentially determined by the measurement accuracy and the sensor placement.

Determining the accuracy of a given sensor configuration is a key prerequisite of optimal sensor placement. 
We develop here a thermal model for a particular machine tool. 
On top of this model we propose two numerical algorithms to evaluate any given thermal sensor configuration with respect to its accuracy. We compute the posterior variances from the posterior covariance matrix with respect to an uncertain initial temperature field. The full matrix is dense and potentially very large, depending on the model size. Thus, we apply a low-rank method to approximate relevant entries, \ie the variances on its diagonal.
We first present a straightforward way to compute this approximation which requires computation of the model sensitivities with with respect to the initial values. Additionally, we present a low-rank tensor method which exploits the underlying system structure.
We compare the efficiency of both algorithms with respect to runtime and memory requirements and discuss their respective advantages with regard to optimal sensor placement problems.

\end{abstract}

% Insert the keywords.

% Insert the Mathematics Subject Classification.
\makeatletter
\ltx@ifpackageloaded{hyperref}{%
\subjclass[2010]{\href{https://mathscinet.ams.org/msc/msc2020.html?t=35K05}{35K05}, \href{https://mathscinet.ams.org/msc/msc2020.html?t=62P30}{62P30}, \href{https://mathscinet.ams.org/msc/msc2020.html?t=68U20}{68U20}, \href{https://mathscinet.ams.org/msc/msc2020.html?t=80A23}{80A23}, \href{https://mathscinet.ams.org/msc/msc2020.html?t=80M10}{80M10}}
}{%
}
\makeatother

% Typeset the opening page.
\maketitle

% Insert the document body.
% International journal for numerical methods in engineering NME
% 
\section{Introduction}
\label{sec:introduction}

Thermal models are nowadays indispensible in the development and operation of the majority of high-precision machinery, for instance, machine tools.
According to \cite{MayrJedrzejewskiUhlmannDonmezKnappHaertigWendtMoriwakiShoreSchmittBrecherWuerzWegener:2012:1}, the thermal error accounts for \SI{75}{\percent} of the total manufacturing error in the final product. 
Although in this paper we focus on applications in the operation of machine tools, the methodology developed can be applied likewise to other classes of transient thermal models.
These models typically contain a large number of uncertain parameters, describing material, process and environmental conditions, such as the heat transfer coefficients between different machine components, ambient temperatures, or the exact heat loss of the engines.
In addition, the initial state, particularly the initial temperature distribution in the machine under consideration, \eg, from its previous operation, are often not known with sufficient accuracy in order to make reliable predictions. 

The majority of the unknowns mentioned above cannot be measured directly but must be inferred from measurements and the underlying thermal model.
The link between the parameters to be estimated during the machine's operation and the available measurements is established via a data assimilation problem \cite{LawStuartZygalakis:2015:1,Peet:2019:1}. 
In our case, the model under consideration is a system of partial differential equations (PDEs) on thermally coupled, complicated geometries. 
These models have to be discretized in space, which lead to large-scale and also stiff state space systems.

In addition to the estimation of the unknown model parameters itself, it is of utmost interest to quantify the impact missing or inexact data has on the estimated parameters and thus on the outcome of the subsequent time-dependent simulation. 
Clearly, carefully positioned measurement devices, collectively referred to as sensors, can help mitigate the effect of inexactness in the measurements.
From personal communication with practitioners, we find that sensors are often placed intuitively, based on an inspection of simulated temperature fields, and in the vicinity of points of interest; see for instance \cite{ThiemKauschingerIhlenfeldt:2018:1,KumarGlaenzelTehelPutz:2021:1}. 
We propose here to follow a systematic approach to assess the impact of measurement errors on the estimation accuracy in dependence on the position of the respective sensor.
To this end, we assume a statistical model with normally distributed, zero-mean measurement errors so that the accuracy is encoded in the estimator's covariance; see, \eg, \cite{CalvettiSomersalo:2007:1}.
To assess the impact on simulation results then requires us to pass this covariance through the forward simulation.

The estimation of these covariances by traditional means involves the setup of large and dense matrices computed from a large number of PDE solutions. 
It therefore requires tremendous computational and memory resources \cite{FlathWilcoxAkcelikHillVanBloemenWaandersGhattas:2011:1,BuiThanhGhattasMartinStadler:2013:1}.
When this is not feasible, a technique to reduce the storage requirements is essential to solve the related posterior covariance problem.
Therefore, we apply a low-rank approximation method to approximate the \textit{data misfit Hessian}, which can then be used to efficiently compute the posterior variances as proposed in \cite{BuiThanhGhattasMartinStadler:2013:1}. 
With this we can greatly reduce storage and computational requirements.
We use tensor train computations to additionally reduce storage requirements during this approximation step.

With regard to the optimal positioning of sensors, we point out that the set of possible sensor locations is typically large.
The evaluation of the respective sensitivity derivatives for an optimal placement strategy is thus computationally expensive and benefits directly from the low-rank approximation as well.

The above considerations apply to any set of unknown parameters.
In this work we consider specifically the quantification of the accuracy of the initial temperature, inferred from solving an estimation problem.
The unknown parameter thus is the infinite-dimensional initial temperature field.
After discretization, the number of parameters agrees with the number of degrees of freedom for the state variable, which is usually quite large.
In turn, the size of the covariance matrix is large as well, which makes the computation even of the variance an expensive task.

For the purpose of this paper, we consider a number of sensors with given locations.
In order to quantify the uncertainty in the parameter estimation problem, the evaluation of sensitivities is required.
Due to the typically moderate number of available sensors, a direct approximation of the required sensitivities is feasible and fast. 
However the storage requirements increase with the number of sensors.
Therefore, in preparation for future work on optimal sensor placement for the problem at hand, we additionally consider a low-rank tensor-train approximation for the solution of the sensitivity equations, making use of the tensor product structure of the underlying system. 
This significantly reduces the storage requirements.
Having both the direct and the low-rank tensor-train approaches at hand allows us to compare them with respect to computational costs and accuracy.

The paper is structured as follows. 
In \cref{sec:the_thermo_elastic_problem} we introduce the PDE system of heat equations to describe the evolution of the temperature fields in a machine tool. 
This section also derives the corresponding state space system by discretization in space, which forms the basis for the direct covariance estimation. 
Statistical background material, particularly concerning the posterior and prior covariance, is collected in \cref{sec:posterior_covariance_matrix}. 
In \cref{sec:computation} we apply the low-rank tensor-train approximation method to the posterior covariance estimation problem. 
This method makes use of the tensor structure of the parameter-to-observable map and does not require storage of the sensitivities with respect to all parameters.
Finally we show in \cref{sec:numerical_results} some numerical results with an emphasis on the comparison of both approaches with regard to time and storage.

\section{Coupled Thermal Models}
\label{sec:the_thermo_elastic_problem}

In this section we describe a class of coupled thermal models which are typical for the description of the temperature state of machine tools.
Each model component is a transient heat equation, \ie, a time-dependent partial differential equation (PDE).
The coupling originates in the heat exchange between different parts of the machine.
This set of PDEs contain material parameters as well as conditions depending on the manufacturing process and the machine's environment.
Each equation has to be discretized in space and time to be solved on a computer.
For concreteness, we now illustrate this on a representative and challenging model problem.

\paragraph{PDE Model of the Auerbach Column}

The \enquote{Auerbach ACW 630} machine is one of the demonstrator objects studied in the collaborative research center \enquote{Thermo-energetic design of machine tools} (SFB/Transregio~96).
The model comprises the machine column, the machine bed and the main spindle; see for instance \cite{NaumannGlaenzelPutz:2020:1}. 

\begin{figure}[htpb]
	\begin{subfigure}[t]{0.45\textwidth}
		\includegraphics[width=0.7\textwidth]{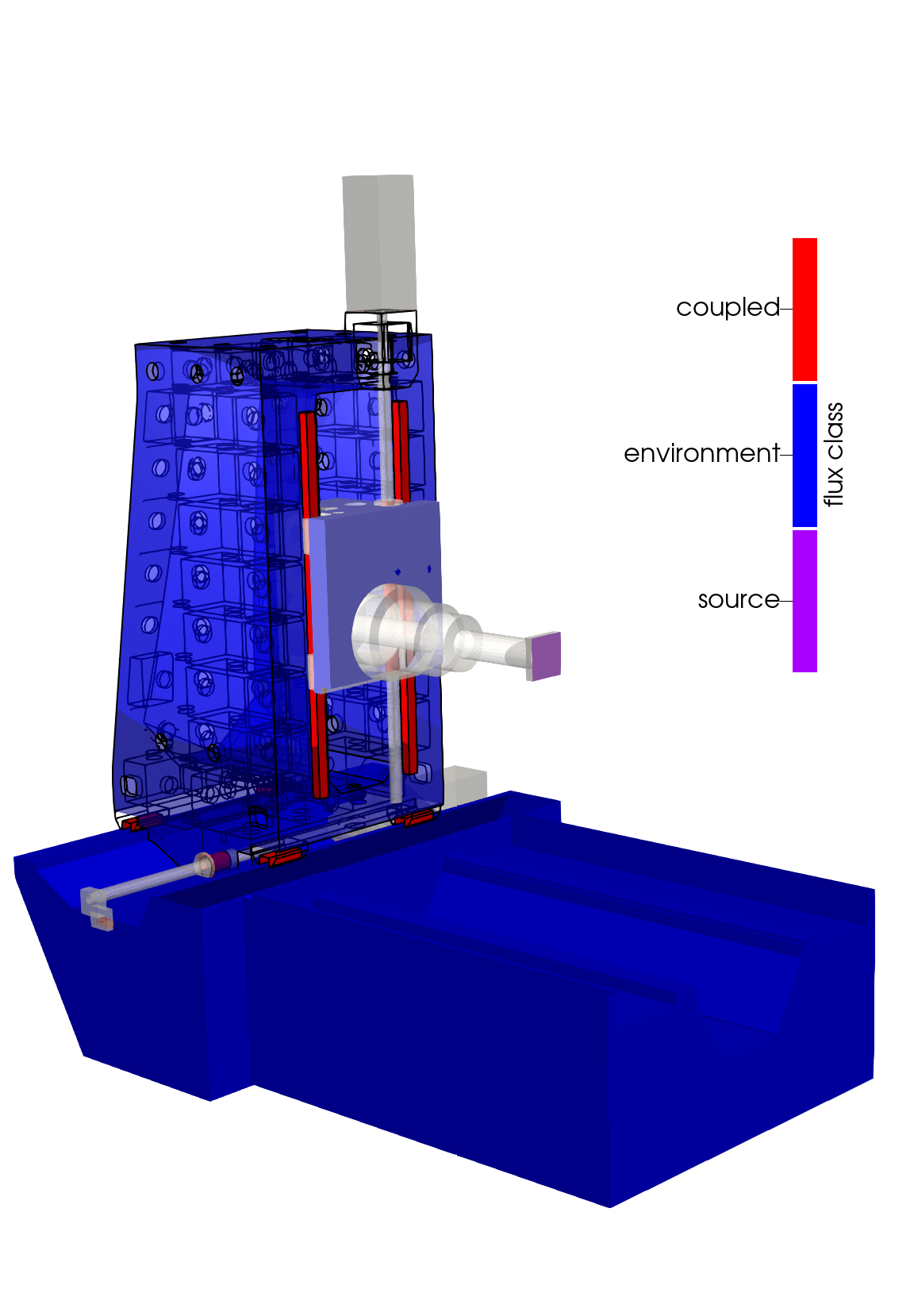}
		\subcaption{The parts and thermal surfaces of the Auerbach machine. The surface consists of three different parts: $\Gamma_\textup{c}$ (red) are coupling surface between machine parts, $\Gamma_\textup{env}$ (blue) is the surface in contact with the environment and $\Gamma_\textup{src}$ (violet) describes the location of heat sources.}%
		\label{fig:ACW:Sources}
	\end{subfigure}
	\hfill
	\begin{subfigure}[t]{0.45\textwidth}
		\includegraphics[width=0.7\textwidth]{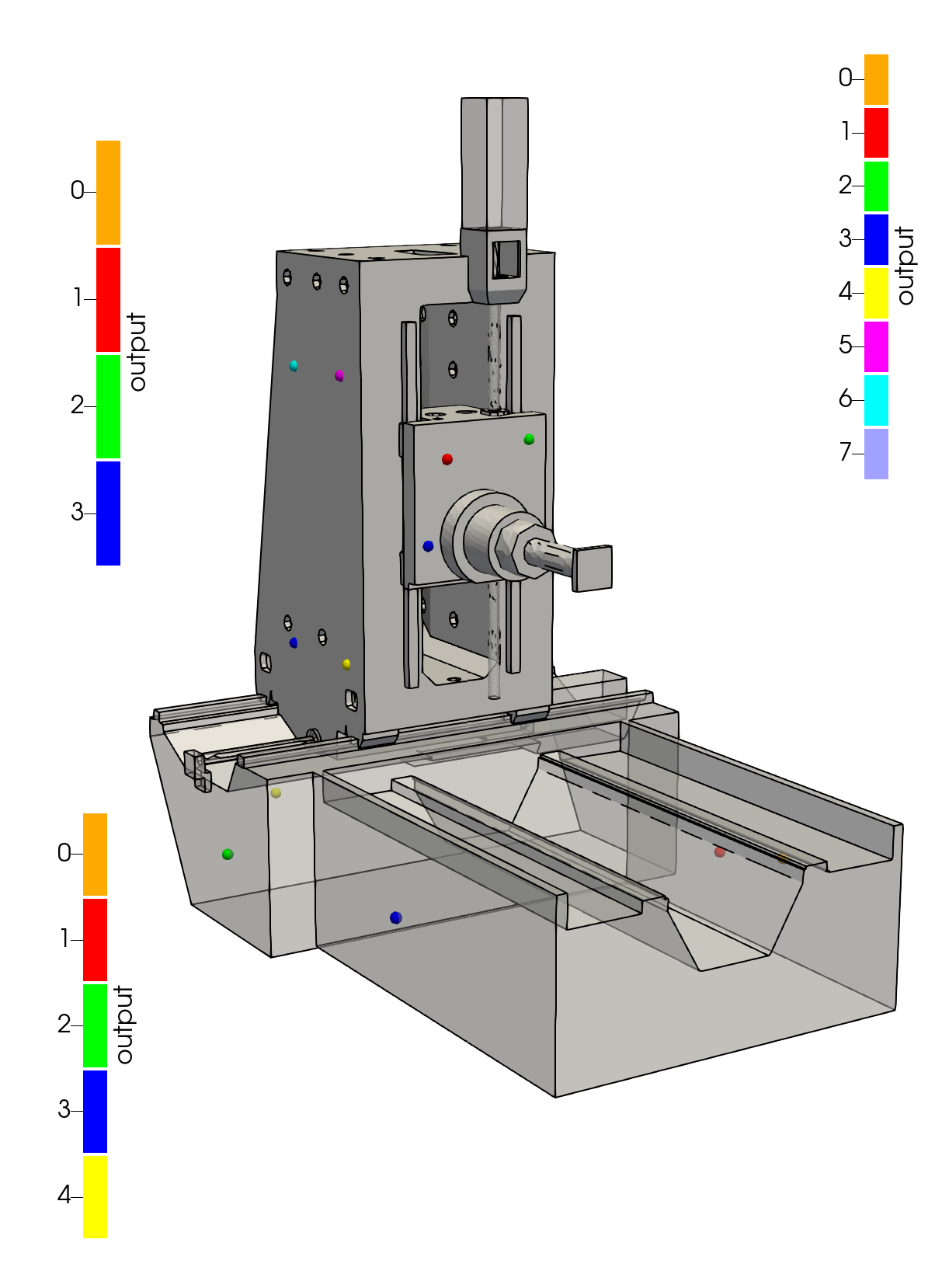}
		\subcaption{In total there are 17~temperature measurement points, representing the outputs $y$ of the model. Five of them are mounted on the machine bed, eight on the machine column, and four on the spindle sledge. Note that not all of them are visible.}
		\label{fig:ACW:Outputs}
	\end{subfigure}
	\caption{Machine geometry showing the different parts of the boundary (left figure), and the measurement points (right figure).}
	\label{fig:ACW}
\end{figure}

The machine is displayed in \cref{fig:ACW}.
Its three parts, the machine column, the bed and the sledge including the spindle, consist of different materials with different densities~$\rho$, thermal capacities~$C_p$ and heat conduction parameters~$\lambda$.
We therefore model each part separately and couple the components by thermal interface conditions.
The volume occupied by part~$i$ is denoted by $\Omega^{(i)}$.
On $\Omega^{(i)}$, the evolution of the temperature field is described by the heat equation
\begin{subequations}
	\label{eq:WLG:ACW}
	\begin{alignat}{4}
		\label{eq:WLG:ACW:volume} 
		\rho^{(i)} \, {C_p}^{(i)} \partial_t T^{(i)} - \div (\lambda^{(i)} \nabla T^{(i)}) 
		&
		= 
		0 
		&
		&
		\quad
		\text{in }
		\Omega^{(i)}
		,
		\\
		\label{eq:WLG:ACW:exchange} 
		n\cdot \lambda^{(i)} \nabla T^{(i)} 
		&
		= 
		q_\textup{c} 
		&
		& 
		\quad
		\text{on } 
		\Gamma_\textup{c} 
		,
		\\
		\label{eq:WLG:ACW:environment} 
		n\cdot \lambda^{(i)} \nabla T^{(i)} 
		&
		= 
		q_\textup{env} 
		& 
		&
		\quad
		\text{on } 
		\Gamma_\textup{env} 
		,
		\\
		\label{eq:WLG:ACW:source} 
		n\cdot\lambda^{(i)} \nabla T^{(i)} 
		&
		= 
		Q_\textup{src}(t) 
		&
		&
		\quad
		\text{on } 
		\Gamma_\textup{src}
		.
	\end{alignat}
\end{subequations}
The first equation \eqref{eq:WLG:ACW:volume} models the temperature evolution inside each machine component.
\Cref{eq:WLG:ACW:exchange} describes the heat exchange with neighboring machine parts at the internal interface $\Gamma_\textup{c}$.
The heat flux into the environment on the boundary $\Gamma_\textup{env}$ is modeled according to \eqref{eq:WLG:ACW:environment}.
The symbol~$n$ denotes the outer normal vector as seen from domain~$\Omega^{(i)}$.
Finally, $\Gamma_\textup{src}$ is the part of the boundary where the model experiences an inflow of thermal energy due to time-dependent external sources $Q_\textup{src}(t)$.

The heat flux $q_\textup{c}$ from machine component~$i$ into component~$j$ is given by
\begin{equation}
	q_\textup{c} 
	=
  \alpha^{(i,j)} \, (T^{(j)} - T^{(i)})
\end{equation}
with a constant exchange coefficient satisfying $\alpha^{(i,j)} = \alpha^{(j,i)}$. 
The heat flux into the environment on $\Gamma_\textup{env}$ is modeled through the term
\begin{equation}
	q_\textup{env} 
	= 
  \alpha_\textup{env} (T^{(j)} - T_\textup{env})
\end{equation}
with heat exchange coefficient $\alpha_\textup{env}$ and constant environmental temperature $T_\textup{env}$.
At the front part of the machine, colored in purple in \cref{fig:ACW:Sources}, a heat source is located, which takes into account thermal energy produced by the respective machining process.
The corresponding part of the boundary, $\Gamma_\textup{src}$, also acts as the tool's point of interest (TCP).
In \cref{fig:ACW:Outputs} we show the location of 17~temperature sensors across the machine's surface.
These sensors are modeled to measure the surface temperature at the point where they are located.

\paragraph{Derivation of the Parameter-to-Observable Map}%
%\label{sub:the_parameter_to_observable_map}

The PDE system \eqref{eq:WLG:ACW} has to be discretized in space. 
We use a standard discretization of the temperature field by linear finite elements \cite[Chapter~8]{ZienkiewiczTaylor:2000:1}.
% \begin{equation*}
%   T^{(i)}(t,x) 
%   = 
%   \sum_{k=1}^{n_x} x^{(i)}_k(t) \, \varphi^{(i)}_k(x)
% \end{equation*}
% of the temperature in the machine component~$i$.
The set of basis functions $\paren[big]\{\}{\varphi^{(i)}_k}_k$ serves as the basis for the spatial discretization of the temperature field in machine part~$\Omega^{(i)}$ and, at the same time, for the  discretization of the unknown initial temperature in that part.
Hence it is natural to define the mass matrices
\begin{equation}
	\label{eq:mass_matrix}
	M^{(i)}_{kl} 
	= 
	\int_{\Omega^{(i)}} \varphi^{(i)}_k \, \varphi^{(i)}_l \d x
	\quad
	\text{and}
	\quad
	\phys{M}^{(i)}_{kl}
	=
	\rho^{(i)} \, {C_p}^{(i)} M^{(i)}_{kl}
	.
\end{equation}
The matrix $M^{(i)}$ serves as the inner product for the parameter space $\R_{M^{(i)}}^{n^{(i)}_x}$ on machine part $\Omega^{(i)}$; see \cite{BuiThanhGhattasMartinStadler:2013:1}.
The total parameter space $\R_M^{n_x}$ is the product space $\R_{M^{(1)}}^{n^{(1)}_x} \times \R_{M^{(2)}}^{n^{(2)}_x} \dots$ with total dimension $n_x = \sum_i n^{(i)}_x$ and the block-diagonal mass matrix $M$, which is comprised of the mass matrices $M^{(i)}$ of all subspaces. 
We remark that the mass matrix represents the $L^2$-inner product in the discrete temperature space.
This is one out of several viable choices, and we anticipate that the inner product carries over to the probability densities and covariance matrices, as will be shown in \cref{sec:posterior_covariance_matrix}.

The discretization of the system \eqref{eq:WLG:ACW} in space, jointly for all three parts of the machine, leads to the affine block input-output ODE system
\begin{subequations}
	\label{eq:DGL}
	\begin{align}
		\phys{M} \dot{x}(t) 
		& 
		= - \phys{K} x(t) + N u(t)
		, 
		\label{eq:DGL1}
		\\
		y(t) 
		&
		= 
		C x(t)
		, 
		\label{eq:DGL2} 
		\\
		x(0) 
		&
		= 
		x_0
		, 
		\label{eq:DGL3}
	\end{align}
\end{subequations}
with $n_x$ spatial degrees of freedom collected in the vector-valued function $x(t)$, an $m$-dimensional input function $u(t)$ combining the heat flux into the system from process sources~$Q_\textup{src}$ and due to the environmental temperature~$T_\textup{env}$, and a $p$-dimensional output function $y(t)$.
As usual, $\dot{x}(t)$ denotes the time derivative of $x(t)$. 
The arising matrices are derived from the finite element discretization of \eqref{eq:WLG:ACW}. 
The symmetric matrices $\phys{M}$ and $\phys{K} \in \R^{n_x \times n_x}$ represent the block-diagonal heat capacity and diffusion matrices including the coupling boundary conditions \eqref{eq:WLG:ACW:exchange}.
The matrix $N \in \R^{n_x \times m}$ describes how each input acts as a distributed heat source.
Finally, the rectangular matrix $C \in \R^{n_y \times n_x}$ transforms the thermal degrees of freedom into the output quantities of interest, \ie, the temperature measurements at the $n_y = 17$~sensor locations.

The initial value problem (IVP) \eqref{eq:DGL} describes the evolution of the temperature state~$x$, and in turn of the output temperatures~$y$ over time, starting from the initial temperature field~$x_0$. 
As mentioned previously, $x_0$ is the uncertain parameter vector.
In order to emphasize the dependence of $y$ on $x_0$, we use the notation~$y(t;x_0)$.

We refer to the outputs~$y$ of the model \eqref{eq:DGL} as the observables. 
The main object of interest for the uncertainty quantification is the parameter-to-observable map~$f$, mapping the initial values~$x_0$ to the outputs~$y(t_s)$ at certain equidistant measurement times~$t_s$, \ie
\begin{equation}
	\label{eq:parameter_to_observables:continous}
	f \colon \R^{n_x}_\phys{M} \to \R^{(n_t + 1) \cdot n_y}
	,
	\quad
	 f(x_0) 
	 = 
	 \begin{pmatrix}
     y(t_0; x_0) 
		 \\
		 y(t_1; x_0) 
		 \\
		 \vdots 
		 \\
		 y(t_{n_t}; x_0)
	 \end{pmatrix}
	 = 
   \sum_{s=0}^{n_t} e_s \otimes y(t_s; x_0)
	 ,
\end{equation}
where $e_s$ is the $s$-th standard unit vector in $\R^{n_t+1}$.
The symbol $\otimes$ denotes the Kronecker product.
In the setting at hand, the parameter-to-observable map is affine-linear.
We denote its linear part, which is obtained when setting $u \equiv 0$ in \eqref{eq:DGL}, by $F \colon \R^{n_x}_M \to \R^{(n_t + 1) \cdot n_y}$.
The (Hilbert space) adjoint of $F$ in the parameter and observation spaces $\R_M^{n_x}$ and $\R^{(n_t + 1) \cdot n_y}$ can be expressed as
\begin{equation}
	\label{eq:adjoint_of_linearized_parameter-to-observable_map}
	F^\natural 
	= 
	M^{-1} F^\transp
	.
\end{equation}

At this point, we have established the basis to quantify the uncertainty in the initial values on the outputs. 
We emphasize that the evaluation of the sensitivities~$F$ of the parameter-to-observable map~$f$ for many points in time and a large number of outputs $n_y \approx n_x$ is not viable in terms of storage.

\section{Posterior Covariance Matrix}
\label{sec:posterior_covariance_matrix}

We recall that it is our goal to infer the initial temperature~$x_0$ from observations of the output~$y(t_s)$ at certain measurement times~$t_s$ for the machine tool and its model \eqref{eq:DGL} described in \cref{sec:the_thermo_elastic_problem}; see \cref{fig:ACW}.
Moreover, we seek to quantify the uncertainty in this estimate.
To this end, we use a Bayesian approach.
In this setting, our prior belief in the probability distribution of the initial state gets informed by measurements to yield a posterior probability distribution.

Our derivation of the background is standard and more details can be found, for instance, in \cite{BuiThanhGhattasMartinStadler:2013:1,Alexanderian:2020:1}.
We temporarily set aside the concrete setting of \cref{sec:the_thermo_elastic_problem} and adopt a general notation.
In this general setting, we have parameters~$P$ of dimension~$n_p$, combined model and observation errors~$E$ and observed outputs~$Y$ of dimension~$n_y$.
All of these are continuous, vector-valued random variables.
As is customary, their realizations are denoted by lower-case letters.
The above are coupled by a (generally nonlinear) forward model
\begin{equation}
	\label{eq:general_forward_model}
	Y 
	= 
	g(P, E)
	.
\end{equation}
We define the following probability density functions (pdfs):
\begin{itemize}
	\item 
		$\pinoise \colon \R^{n_y} \to \R$, representing the model error and observation noise,
	\item 
		$\piprior \colon \R^{n_p}_M \to \R$, describing the prior information about the model parameters~$P$,
\end{itemize}
and the likelihood function $\pi(y | p)$, describing the probability of the observable~$y$, given the model parameters~$p$.
Note that the parameter space $\R^{n_p}_M$ is equipped with an inner product $\inner{u}{v}_M = u^\transp M \, v$ represented by a symmetric, positive definite matrix~$M$, as is the case for our model; see \eqref{eq:mass_matrix}.

Exploiting Bayes' theorem, we can combine the prior probability density and the likelihood function to form the posterior probability density $\pipost \colon \R^{n_p}_M \to \R$ as follows,
\begin{equation}
	\pipost(p) 
	= 
	\pi(p|y) 
	= 
	\frac{\piprior(p) \, \pi(y|p)}{\pi(y)} 
	\propto 
	\piprior(p) \, \pi(y|p)
	.
\end{equation}
We can further simplify this by assuming additive noise, in which case \eqref{eq:general_forward_model} is replaced by
\begin{equation}
 \label{eq:parameter_to_observable} 
	Y 
	= 
	f(P) + E
\end{equation}
with the noise-free parameter-to-observable map~$f$.
We can then write $\pinoise(e) = \pinoise(y - f(p))$ and with this we get
\begin{equation} 
	\label{eq:post_density}
	\pipost(p)
	\propto 
	\piprior(p) \, \pinoise(y - f(p))
	.
\end{equation}
As is often the case, we assume that both the prior and noise probability densities are Gaussian, \ie, 
\begin{align}
	\label{eq:prior_Gaussian}
	\pi_\textup{prior}(p) 
	&
	\propto 
	\exp \paren[Big](){-\frac{1}{2} \norm{p - \bar{p}}_{\gPrior^{-1}}^2} 
	,
	\\
	\label{eq:noise_Gaussian}
	\pinoise(e) 
	&
	\propto \exp \paren[Big](){- \frac{1}{2} \norm{e - \bar{e}}_{\gNoise^{-1}}^2}
	.
\end{align}
Here, $\bar{p}_\textup{prior} \in \R^{n_p}_M$ denotes the mean of the model parameter's prior pdf and $\bar{e} \in \R^{n_y}$ is the mean of the noise pdf.
Moreover, $\gPrior \in \R^{{n_p} \times {n_p}}$ is the prior covariance, which is positive definite and symmetric.
Moreover, $\gNoise \in \R^{n_y \times n_y}$ is the symmetric positive definite covariance matrix of the noise.

With the Gaussian densities \eqref{eq:prior_Gaussian} and \eqref{eq:noise_Gaussian} we can write the posterior probability density \eqref{eq:post_density} as
\begin{align*}
	\pipost(p) 
	&
	\propto \exp \paren[auto](){-\frac{1}{2} \paren[big][]{\norm{p-\bar{p}_\textup{prior}}_{\gPrior^{-1}}^2 + \norm{e-\bar{e}}_{\gNoise^{-1}}^2}} 
	\\
	&
	= 
	C \exp \paren[auto](){-\frac{1}{2} \paren[big][]{\norm{p-\bar{p}_\textup{prior}}_{\gPrior^{-1}}^2 + \norm{y-f(p)-\bar{e}}_{\gNoise^{-1}}^2}}
\end{align*}
with a normalization constant $C$ which is independent of $p$.
That said, the posterior probability density $\pipost(p)$ is not necessarily Gaussian as the function $f(p)$ may be nonlinear. 
Recall that in our setting \eqref{eq:DGL}--\eqref{eq:parameter_to_observables:continous}, however, we have an affine map of the form $f(p) = F \, p + f_0$, which takes the initial state~$x_0$ to the observed outputs.
In this case we obtain that $\pipost(p)$ is Gaussian with mean $\bar{p}_\textup{post}$ given by the maximum a~posteriori probability (MAP) point
\begin{align}
	\bar{p}_\textup{post} 
	&
	= 
	\argmax_{p \in \R^{n_p}_M} \pipost(p) 
	\notag
	\\
	&
  = 
	\argmin_{p \in \R^{n_p}_M} \frac{1}{2} \paren[big][]{\norm{p-\bar{p}_\textup{prior}}_{\gPrior^{-1}}^2 + \norm{y - Fp - f_0 - \bar{e}}_{\gNoise^{-1}}^2}
	. 
	\label{eq:map}
\end{align}
The quantity we are primarily interested in is the covariance of the posterior pdf, $\gPost$. 
It is given by the inverse of the Hessian of the negative logarithm of $\pipost$, \ie, 
\begin{equation}
	\gPost
  = 
	\paren[auto](){{F^{\transp}} \gNoise^{-1}F + \gPrior^{-1}}^{-1}
	.
	\label{eq:post_covar}
\end{equation}
This quantity gives us a measure for the posterior covariance in terms of the noise and prior covariances and depending on the linear part of the parameter-to-observable map, which in turn depends on the outputs we choose.

The noise covariance $\gNoise$ will be specified along with the numerical examples in \cref{sec:numerical_results}.
For the prior, we employ a Laplacian-like covariance structure which is in accordance with the infinite-dimensional nature of the initial state.
It will be discussed in more detail in \cref{sec:prior_selection}.

The posterior covariance \eqref{eq:post_covar} itself is usually infeasible to work with, as, in our setting, it is a large and dense matrix of dimension~$\gPost \in \R^{n_x \times n_x}$. 
Even though our particular interest lies only with the posterior variances, \ie, the diagonal entries of $\gPost$ given by
\begin{equation}
  \label{eq:post_variance}
	\text{var}_k 
	= 
  e_k^\transp \gPost \, e_k
\end{equation}
with unit vectors $e_k \in \R^{n_x}$, a storage reduction in the evaluation of $\gPost$ is crucial. 
Fortunately, a low-rank approximation method has emerged, which delivers reliable results while significantly reducing the computational requirements.
We describe it in the following subsection.

\subsection{\texorpdfstring{Low-Rank Approximation Method for $\gPost$}{Low-Rank Approximation Method for the Posterior}}
\label{subsection:Low-rank_approximation}

The posterior covariance \eqref{eq:post_covar} can be approximated efficiently as explained in the following.
For convenience, we rewrite \eqref{eq:post_covar} as
\begin{equation}
  \gPost 
	= 
	\paren[auto](){\gPrior F^\transp \gNoise^{-1} F + \id}^{-1} \gPrior,
\end{equation}
where $\hMisfit \coloneqq F^\transp \gNoise^{-1} F$ is the Hessian of the model part in \eqref{eq:map}.
In addition we call $\gPrior \hMisfit$ the preconditioned misfit Hessian.
We aim for a low-rank approximation based on the spectral decomposition of the preconditioned misfit Hessian of the form 
\begin{equation}
	\label{eq:HmisApprox} 
	\gPrior \hMisfit 
	=
	V \Lambda V^\dag
	\approx 
	V_r \Lambda_r V_r^\dag
	.
\end{equation}
The columns of $V$ contain the eigenvectors, $V^\dag = V^\transp \gPrior^{-1}$ is the Hilbert-space adjoint of $V$, and the diagonal matrix $\Lambda$ contains the corresponding eigenvalues of the eigenvalue problem
\begin{equation*}
  \gPrior \hMisfit \, v 
	= 
	\lambda \, v
	,
\end{equation*}
which is equivalent to the generalized eigenvalue problem
\begin{equation}
 \label{eq:hMis:genEig} 
 \hMisfit \, v
 = 
 \lambda \, \gPrior^{-1} v
 .
\end{equation}
The misfit Hessian is symmetric and the prior covariance is symmetric and positive definite, and thus the eigenvectors~$v$ are orthogonal with respect to $\inner{\cdot}{\cdot}_{\gPrior^{-1}}$, \ie, $V^\transp \gPrior^{-1} V = \id$, and thus we also have $V \, V^\transp \gPrior^{-1} = V \, V^\dag = \id$.
The low-rank approximation of rank~$r$ is built using the leading eigenpairs.
We can compute an approximate generalized spectral decomposition \eqref{eq:hMis:genEig} of $\hMisfit$, \eg, via the Arnoldi method with respect to the $\gPrior^{-1}$-inner product.
The rank~$r$ can be adjusted depending on the available computing time and memory resources.
  
Using this low-rank decomposition, we can calculate an approximate inverse of $\gPrior\hMisfit + \id$ using the orthogonality $V V^\dag = \id$, which yields
\begin{equation}
	\label{eq:gammaPostInverse} 
  \paren[big](){\gPrior\hMisfit + \id}^{-1}
   = 
	 \id - V \widetilde{\Lambda} V^\dag 
  \approx 
	\id - V_r \widetilde{\Lambda}_r V_r^\dag
	,
\end{equation}
with $\widetilde{\Lambda} = \diag\paren[auto](){\frac{\lambda_1}{\lambda_1 + 1}, \ldots, \frac{\lambda_{(n_t+1)n_y}}{\lambda_{(n_t+1)n_y} + 1}}$.
Finally, the posterior covariance matrix $\gPost$ and its diagonal entries can be evaluated according to
\begin{equation}
	\label{eq:evaluation_of_diagonal_entries_of_posterior_covariance}
	e_k^\transp \gPost \, e_k
	\approx
  e_k^\transp \paren[big](){\id - V_r \widetilde{\Lambda}_r V_r^\dag } \gPrior \, e_k
	=
  e_k^\transp \gPrior \, e_k - e_k^\transp V_r \widetilde{\Lambda}_r V_r^\transp e_k
	.
\end{equation}

Let us relate \eqref{eq:hMis:genEig}--\eqref{eq:evaluation_of_diagonal_entries_of_posterior_covariance} to the representation given in \cite[Sections~3.7, 5.2 and 5.3]{BuiThanhGhattasMartinStadler:2013:1}.
Their prior $\gPriorO$ is related to ours by $\gPrior = \gPriorO M^{-1}$.
From the appearance of the inner product~$M$ in the parameter space, it becomes apparent that they represent the prior as an endormorphism in the parameter space $\R^{n_p}_M$, while we represent the same as a linear map from the dual space $(\R^{n_p}_M)^*$ into the parameter space.
Consequently, while our prior representation is positive definite and self-adjoint, theirs is positive definite and Hilbert space self-adjoint, \ie
\begin{equation*}
	M^{-1} \gPriorO^\transp
	=
	\gPriorO M^{-1} 
	.
\end{equation*}
The authors in \cite{BuiThanhGhattasMartinStadler:2013:1} further assumed that the square root, $\gPriorO^{1/2}$, is accessible, \ie, we have
\begin{equation}
	\label{eq:rootO}
	\gPriorO 
	= 
	\gPriorO^{1/2} \gPriorO^{1/2}
	,
\end{equation}
where $\gPriorO^{1/2}$ is likewise positive definite and Hilbert space self-adjoint, \ie $M^{-1} \gPriorO^{\transp/2} = \gPriorO^{1/2} M^{-1}$ holds.

We proceed to provide an alternative expression for the generalized eigenvalue problem \eqref{eq:hMis:genEig} and thus for the representation \eqref{eq:evaluation_of_diagonal_entries_of_posterior_covariance}.
We begin with \eqref{eq:post_covar}, which we now write as
\begin{align*}
	\gPost 
	&
	= 
	\paren[auto](){M^{-1}\hMisfit + \gPriorO^{-1}}^{-1}M^{-1} 
	\\
	&
	=
	\gPriorO^{1/2}\paren[auto](){\gPriorO^{1/2} M^{-1} \hMisfit \gPriorO^{1/2} + \id}^{-1} \gPriorO^{1/2}M^{-1}
	\\
	&
	=
	\gPriorO^{1/2}\paren[auto](){M^{-1} \gPriorO^{T/2} \hMisfit \gPriorO^{1/2} + \id}^{-1} \gPriorO^{1/2}M^{-1}
	,
\end{align*}
where we applied the Hilbert space self-adjointness of $\gPriorO^{1/2}$ in the last step.
  From here we repeat the steps leading from \cref{eq:HmisApprox} to \cref{eq:evaluation_of_diagonal_entries_of_posterior_covariance}.
	This time, we obtain the generalized eigenvalue problem
	\begin{equation}
		\label{eq:hMis:genEigO}
    \gPriorO^{T/2} \hMisfit \gPriorO^{1/2}\hat{v} 
		= 
		\hat{\lambda} \, M \, \hat{v}
  \end{equation}
	and spectral decomposition
  \begin{equation}
		\label{eq:hMis:genEigO:decomposition}
    M^{-1} \gPriorO^{T/2} \hMisfit \gPriorO^{1/2} 
		= 
		\hat{V} \Lambda \hat{V}^\#
		.
  \end{equation}
	It is easy to see that the eigenvalues are the same as in \eqref{eq:hMis:genEig} and the eigenvectors are related through $\hat{v} = \gPriorO^{-1/2} v$.
		Moreover, we have $\hat{V}^\# = \hat{V}^\transp M$, and the eigenvectors $\hat{v}$ are orthogonal with respect to $M$.
	As above, we can use the leading $r$~eigenpairs to obtain the approximation
  \begin{equation*}
    \paren[auto](){M^{-1}\gPriorO^{T/2}\hMisfit\gPriorO^{1/2} + \id}^{-1} 
		\approx 
		\id - \hat{V}_r \widetilde{\Lambda}_r \hat{V}^\#_r
		.
  \end{equation*}
  Thus the expression for the variance, analogous to \eqref{eq:evaluation_of_diagonal_entries_of_posterior_covariance}, reads 
	\begin{equation}
		\label{eq:evaluation_of_diagonal_entries_of_posterior_covarianceO}
		e_k^\transp \gPost e_k 
		\approx 
		e_k^\transp \gPrior e_k - e_k^\transp \gPriorO^{1/2} \hat{V}_r \widetilde{\Lambda} \hat{V}_r^\transp \gPriorO^{T/2} e_k
		.
	\end{equation}
	We thus conclude that both representations of the approximate posterior variance differ only \wrt the formulation of the respective generalized eigenvalue problems.
	However, both formulations eventually lead to the same approximation of the posterior variance.

\subsection{Construction of the Prior Covariance}%
\label{sec:prior_selection}

The Bayesian approach requires the user to specify their prior knowledge on the unknown parameter.
This is specifically expressed through the prior covariance operator or matrix $\gPrior$ for the initial temperature. 
Our approach is analogous to \cite{BuiThanhGhattasMartinStadler:2013:1}.
Using the stationary variant of \eqref{eq:WLG:ACW} as a model, we assume a covariance operator for the initial temperature based on the square of an inverse Laplace operator $A = b - a \, \laplace$ with parameters $a, b > 0$.

Next we motivate our choice for these parameters by considering a Dirac-type initial state $T(0,\cdot) = \delta_x$.
In the absence of heat sources, an implicit Euler step 
\begin{equation*}
	\paren[auto](){\rho \, C_p \, \id - \tau \, \lambda \, \Delta} T(\tau, \cdot) 
	= 
	\rho \, C_p T(0, \cdot) 
\end{equation*}
gives us an estimate of the temperature profile after a short time~$\tau$.
In simple geometries, the ratio $\beta = \frac{b}{a}$ controls the decay of the solution away from the location of the Dirac pulse; see, \eg, \cite[section~8.3.2]{Polyanin:2002:1}.
With $\rho$, $C_p$ and $\lambda$ given material parameters, we prescribe the time constant $\tau = \SI{30}{\minute}$ and thus we have a computable constant
\begin{equation*}
	\beta
	=
	\frac{b}{a}
	=
	\frac{\rho \, C_p}{\tau \, \lambda}
	.
\end{equation*}
% We choose the value of $b$ (and deduce $a$) further below.
Moreover, we choose $a$ so that the average prior variance over the finite element mesh nodes is prescribed and is equal to $\SI{3}{\kelvin\squared}$.
To this end, we discretize $A = b - a \, \laplace$ with linear finite elements on the mesh, \ie, we have 
\begin{equation}
	\label{eq:prior:discretePDE}
  \paren[big][]{A_h}_{k\ell}
	= 
  \int_{\Omega} a \nabla \varphi_k \cdot \nabla \varphi_\ell + b \, \varphi_k \, \varphi_\ell \d x 
	=
  aK + bM
	.
\end{equation}
Here $K$ and $M$ are the stiffness and mass matrices similar to those in \eqref{eq:DGL} but with unit coefficients.
The prior is related to $A_h$ as follows,
\begin{equation}
	\gPrior 
	= 
	A_h^{-1} M A_h^{-1}
	.
\end{equation}
The variance in each finite element node is equal to a diagonal entry of $\gPrior$, and depends on the unknown~$a$ and the chosen constant $\beta$ according to 
\begin{equation*}
	A_h 
	= 
	a \, (K + \beta \, M)
\end{equation*}
We then choose $a$ such that the average prior variance over the finite element mesh nodes is equal to $\SI{3}{\kelvin\squared}$.
The coefficient $b = \beta \, a$ then follows as well.

In our application example concerning the machine tool from \cref{sec:the_thermo_elastic_problem}, which consists of three parts, we apply this reasoning individually on each part separately and thus obtain parameters $a^{(i)}$ and $b^{(i)}$.
We refer the reader to \cref{sec:numerical_results}, \cref{fig:priorVariance:physicalParameters} for a visualization of the prior's pointwise variance as well as for a comparison of the prior and posterior variance fields.

\begin{remark*}
	The above made construction ensures that the matrix square root \eqref{eq:rootO} of $\gPriorO$ is explicitly available, namely we have
	\begin{equation*}
		\gPriorO^{1/2} 
		=
		A_h^{-1} M
		.
	\end{equation*}
	This allows us to use the alternative representation \eqref{eq:hMis:genEigO}.
\end{remark*}

\section{Efficient Evaluation of the Linearized Parameter-to-Observable-Map}
\label{sec:computation}

We recall that it is our goal to evaluate the posterior nodal variances \eqref{eq:evaluation_of_diagonal_entries_of_posterior_covariance} of the initial temperature field according to \eqref{eq:evaluation_of_diagonal_entries_of_posterior_covariance}. 
We also recall that this requires the determination of the $r$~leading eigenpairs of the prior-preconditioned misfit Hessian $\gPrior F^\transp \gNoise^{-1} F$; see \eqref{eq:HmisApprox}.
This is achieved via an Arnoldi iteration in the $\gPrior^{-1}$-inner product, which requires the repeated evaluation of matrix-vector products with $F^\transp \gNoise^{-1} F$.
This routine is implemented in \arpack, see \cite{LehoucqSorensenYang:1998:1}.
Recall that $F$ is the linear part of the parameter-to-observable map; see the text following \eqref{eq:parameter_to_observables:continous}.

In this work we present two approaches to evaluate these matrix-vector products.
In the first, straightforward approach, we pre-compute the dense matrix~$F \in \R^{(n_t+1)n_y \times n_x}$.
This is only feasible when both the mesh size~$n_x$ and the number of outputs~$n_y$ and time steps~$n_t$ are moderate.
We refer the reader to \cref{sub:direct_approximation} for this direct approach.

An alternative is to replace matrix-vector products with $F$ and with $F^\transp$ using a tensor-train approximation.
This is feasible even for large-scale problems where the direct approach would fail.
We discuss the tensor-train approach in \cref{sub:low_rank_time_discretization}.

\subsection{Direct Approximation}
\label{sub:direct_approximation}

For the application at hand, the linear part $F \in \R^{(n_t+1)n_y \times n_x}$ of the parameter-to-observable map is typically a matrix with only a small number of rows, while the column dimension represents the discretization of the initial value and is large.
Therefore, we use an adjoint approach to evaluate $F$.
For this approach we derive a backwards-in-time adjoint system, where each time step pertains to one block-row of dimension $n_y = 17$ in $F$.
Notice that since the forward problem \eqref{eq:DGL} is affine, the linearization $F$ does not depend on the state itself.

In view of \eqref{eq:DGL}, the continuous adjoint problem to be solved is the matrix-valued final value problem (FVP)
\begin{subequations}
	\label{eq:sensitivity}
	\begin{align}
		\label{eq:sensitivity:factors}
		\phys{M}^\transp\dot{S}(t) 
		&
		=
		\phys{K} S(t) 
		\\
		\phys{M}^\transp S(T) 
		&
		= 
		-C^\transp
	\end{align}
\end{subequations}
for the sensitivity factors $S(t) \in \R^{n_x \times n_y}$.
With $S$ known, the sensitivity of the outputs with respect to the initial values are given by
\begin{equation}
	\label{eq:sensitivity:initial}
	\frac{\d y(t_s)}{\d x_0}
	=
	- S(T - t_s)^\transp \phys{M}
	.
\end{equation}
Notice that due to the autonomous structure of the forward system governing $F$, it is enough to solve only a single adjoint system (as opposed to one adjoint system per observation time). 
We refer the reader for instance to \cite[Sect.~2]{HerzogRiedelUcinski:2018:1} for a similar approach in the context of thermo-elastic problems.
In analogy to the definition of the parameter-to-observable-map, we stack the sensitivities \eqref{eq:sensitivity:initial} to obtain 
\begin{equation}
	\label{eq:derivative:continous}
	F 
	= 
	-
	\begin{bmatrix}
		S(T-0)^\transp \phys{M}
		\\
		\vdots
		\\
		S(T-T)^\transp \phys{M}
	\end{bmatrix}
	=
	- \sum_{s = 0}^{n_t} e_s \otimes [S(T -t_s)^\transp \phys{M}]
	,
\end{equation}
where $e_s$ is the $s$-th standard unit vector in $\R^{n_t+1}$.
In our numerical implementation, we solve \eqref{eq:sensitivity} using the implicit Euler method with $n_t$~equidistant time steps of length~$\tau$ identical to the measurement times.
This amounts to 
\begin{equation}
	\label{eq:fully_discrete_adjoint_problem}
	\frac{\phys{M}^\transp S(t_s) - \phys{M}^\transp S(t_{s-1})}{\tau} 
	= 
	- \phys{K} S(t_{s-1})
	, 
	\quad 
	\text{ for } 
	s = n_t, \ldots, 1
\end{equation}
with $t_s = s \, \tau$, $s = 0, \ldots, n_t$, and terminal values~$S(T) = -C^\transp$.

\subsection{Tensor Train Approximation}%
\label{sub:low_rank_time_discretization}

The tensor-train approximation requires us to take a closer look at the structure of the affine forward system \eqref{eq:DGL}.
Also here we use an implicit Euler scheme for the time discretization with $n_t$ equidistant time steps of length~$\tau$ to discretize the time interval $t \in [0,T]$ into steps $t_s = s \, \tau$ for $s = 0,\ldots, n_t$. 
With this we obtain the fully discrete system
\begin{subequations}
	\label{eq:fully_discrete_forward_problem}
	\begin{align}
		\frac{\phys{M} x_s - \phys{M} x_{s-1}}{\tau} 
		&
		= 
		- \phys{K} x_s + N u_s
		, 
		\quad 
		\text{ for } i=1,\ldots, n_t 
		\\
		y_s 
		&
		= 
		C x_s,
	\end{align}
\end{subequations}
with initial values~$x_0$.
Here, $x_s$, $y_s$ and $u_s$ are the discrete values of the functions $x(t)$, $y(t)$ and $u(t)$ at time~$t_s$. 
We denote the matrix collecting these vectors at all time steps by 
\begin{equation*}
	\bx 
	= 
	\begin{bmatrix} 
		x_0
		\\
		x_1
		\\
		\vdots
		\\
		x_{n_t} 
	\end{bmatrix}
	\in \R^{n_x (n_t+1)}
\end{equation*}
and similarly for $\by$ and $\bu$.

In order to motivate a low-rank tensor approximation scheme for \eqref{eq:fully_discrete_forward_problem}, we formally write this system as one large linear system of equations,
\begin{subequations}
	\begin{align}
		\MoveEqLeft
		\underbrace{%
			\begin{bmatrix} 
				\tau \phys{K} + \phys{M} 
				& 
				& 
				\\
				- \phys{M} 
				& 
				\tau \phys{K} + \phys{M} 
				& 
				\\ 
				& 
				\hspace{-3em} \ddots 
				& 
				\hspace{-3em} \ddots 
				\\ 
				& 
				-\phys{M} 
				& 
				\tau \phys{K} + \phys{M} 
			\end{bmatrix}
		}_{\tens{K}} 
		\begin{bmatrix} 
			x_0 
			\\
			x_1 
			\\ 
			\vdots 
			\\ 
			x_{n_t} 
		\end{bmatrix} 
		\notag
		\\
		&
		= 
		\underbrace{\begin{bmatrix} N & & \\ & \ddots & \\ & & N \end{bmatrix}}_{\tens{N}} \begin{bmatrix} 0 \\ u_1 \\ \vdots \\ u_{n_t} \end{bmatrix} + \underbrace{\begin{bmatrix} \tau \phys{K} + \phys{M} & & & \\ & 0 & & \\ & & \ddots	& \\ & & & 0 \end{bmatrix}}_{\tens{M}_0} \underbrace{\begin{bmatrix} x_0 \\ 0 \\ \vdots \\ 0 \end{bmatrix}}_{\tens{x}_0}
		, 
		\\
		\begin{bmatrix} y_0 \\ \vdots \\ y_{n_t} \end{bmatrix} 
		&
		= 
		\underbrace{\begin{bmatrix} C & & \\ & \ddots & \\ & & C \end{bmatrix}}_{\tens{C}} \begin{bmatrix} x_0 \\ \vdots \\ x_{n_t} \end{bmatrix}
		.
	\end{align}
\end{subequations} 
Using this formulation we can denote the parameter-to-observable map $x_0 \mapsto \tens{y}$ concisely as
\begin{equation} 
	\label{eq:parameter_to_observables}
	\tens{y} 
	= 
	\tens{C} \tens{K}^{-1} \paren[auto](){\tens{N} \tens{u} + \tens{M}_0 \tens{x}_0}. 
\end{equation}
The matrices involved can be expressed as the following Kronecker products,
\begin{equation}
	\label{eq:kronecker} 
	\begin{alignedat}{3}
			\tens{C} 
			&
			= 
			\id_{n_t} \otimes C
			,
			\quad
			&
			\tens{N} 
			&
			= 
			\id_{n_t} \otimes N
			, 
			\\
			\tens{K} 
			&
			= 
			\tau \, \id_{n_t} \otimes \phys{K} + B \otimes \phys{M}
			,
			\quad
			&
			\tens{M}_0 
			&
			= 
			\id_0 \otimes (\tau \phys{K} + \phys{M})
			, 
		\end{alignedat}
	\end{equation}
where $\id_{n_t}$ is the identity matrix of dimension $(n_t+1) \times (n_t+1)$, $B$ is a bidiagonal matrix with $1$ on the diagonal and $-1$ on its first subdiagonal and $\id_0$ is an almost all-zero matrix with its $(1,1)$ entry equal to $1$.

While the system matrix $\tens{K}$ is prohibitively large, it is full of structure that we will now show can be used to develop efficient numerical schemes. 
Also, $\tens{C}$ and $\tens{M}_0$ have a Kronecker product structure that can more readily be exploited, given that they do not consist of a sum of components. All-at-once systems, such as the one contained in $\tens{K}$, have also been discussed in \cite{StollBreiten:2015:1,HerzogPearsonStoll:2019:1,PearsonStollWathen:2012:1,BenziHaberTaralli:2011:1}.
We recall that we require matrix-vector products with $F^\transp \gNoise^{-1} F$ in order to compute an approximate spectral decomposition as in \eqref{eq:HmisApprox}.
Notice that $F$ is given by \eqref{eq:parameter_to_observables} with $\tens{u} = 0$, \ie,
\begin{equation*}
	F 
	= 
	\tens{C} \tens{K}^{-1} \tens{M}_0
	.
\end{equation*}

In contrast to the direct approach in \cref{sub:direct_approximation} we now exploit the Kronecker product structure \eqref{eq:kronecker}.
We follow \cite{BennerQiuStoll:2018:1} where the authors propose using the tensor-train (TT) format \cite{Oseledets:2011:1,OseledetsTyrtyshnikov:2009:1} to perform these calculations efficiently. 
The TT format is especially convenient for our purpose as the underlying Kronecker product structure \eqref{eq:kronecker} of the parameter-to-observable map \eqref{eq:parameter_to_observables} can be interpreted as a low-rank tensor-train representation. 

We proceed to describe our approach to form approximate matrix-vector products $F v$.
The case for $F^\transp$ is then similar.
Owing to the Kronecker structure of $\tens{C}$ and $\tens{M}_0$, products with these matrices involve only matrix-vector products of size $n_x \times n_x$, as they do not contain any sum of terms and no system with them needs to be solved.

We do require matrix-vector products with $\tens{K}^{-1}$, which we view as the solutions of $\tens{K} \tens{x} = \tens{b}$. 
We now briefly discuss how such systems can be solved within a general tensor framework.
The unknown $\tens{x}$ will be represented in tensor-train format, which allows it to be represented with adjustable rank.
This format permits an efficient iterative solution of $\tens{K} \tens{x} = \tens{b}$, using only tensor contractions in place of matrix-vector products.

In our setting, we use a two-dimensional tensor-train format, distinguishing the temporal from the spatial directions but additional structure such as Kronecker structures in mass and stiffness matrices can easily be incorporated (\cite{BuengerDolgovStoll:2020:1}). 
In our case, $\tens{x}$ will be represented in tensor-train format as
\begin{equation}
	\label{eq:TT_representation}
	\tens{x} 
	= 
	\sum_{\alpha=1}^{r} x_{\alpha}^{(1)} \otimes x_{\alpha}^{(2)}
\end{equation}
with the cores $x_{\alpha}^{(1)} \in \R^{n_t+1}$ and $x_{\alpha}^{(2)}\in \R^{n_x}$, and $r$ is a suitably chosen low-rank parameter.
Such a format is particularly attractive for a large number of time steps and/or large number of output points due to the considerable reduction in storage cost for moderate rank~$r$.

We view $\tens{F} \in \R^{(n_t+1)n_y \times n_x}$ as a tensor approximated in the tensor-train format. 
Obviously, we are only interested in the efficient evaluation of $\tens{F} \tens{v}$ where $\tens{F}$ and $\tens{v}$ are represented in the tensor-train format. 
Within this process the solution with $\tens{K}$ to obtain $\tens{x}$ via the system $\tens{K} \tens{x} = \tens{b}$ is the most challenging part and we discuss this in some detail now.
We review the state-of-the-art solver in general form and point to the literature for detailed descriptions. This can be done efficiently with an energy function minimization such as the \textit{alternating linear scheme} (ALS) \cite{HoltzRohwedderSchneider:2012:1}. 
This approach constructs low-dimensional systems of linear equations for each core which can then be solved with standard numerical methods. 
Here, we derive the method for symmetric $\tens{K}$. 
Solving $\tens{K} \tens{x} = \tens{b}$ then corresponds to the minimization of the energy function,
\begin{equation}
	\min_{\tens{x}} J(\tens{x}) 
	= 
	\norm{\tens{x}_* - \tens{x}}^2_\tens{K} 
	= 
	(\tens{x}, \tens{K} \tens{x}) - 2 (\tens{x}, \tens{b}) + \text{const}
	, 
	\label{eq:energy}
\end{equation}
with the exact solution $\tens{x}_* = \tens{K}^{-1} \tens{b}$. 
We will now compute the solution to our problem $\tens{x}$ in the form \eqref{eq:TT_representation}.
To find a solution for \eqref{eq:energy} we select an initial guess $\tens{x}_0$ and cycle over its TT-cores where we solve a reduced version of \eqref{eq:energy} to improve the current guess. 
The alternating solver proceeds by updating the solution with respect to one of the dimensions freezing the other cores in the process. 
For this we rely on the linearity of the tensor-train format, \ie,
\begin{equation}
  \tens{x} 
	= 
	\tens{x}_{\neq k} x_{=k}
	,
\end{equation}
with $x_{=k}$ a vectorization of core $\tens{x}^{(k)}$ and $\tens{x}_{\neq k}$ the tensor-train where $\tens{x}^{(k)}$ is replaced by a placeholder identity of appropriate size. 
We use this notation to construct reduced problems to sequentially update the cores individually.
With this, the energy function for the reduced problem for the $k$-th core becomes
\begin{equation}
		\label{eq:local:energy}
	\begin{aligned}
		J(\tens{x}) 
		&
		= 
		(\tens{K} \tens{x}_{\neq k} x_{=k}, \tens{x}_{\neq k} x_{=k}) - 2(\tens{b},\tens{x}_{\neq k} x_{=k}) 
		\\
		&
		= 
		(\tens{x}_{\neq k}^* \tens{K} \tens{x}_{\neq k} x_{=k},x_{=k}) - 2(\tens{x}_{\neq k}^* \tens{b},x_{=k})
		.
	\end{aligned}
\end{equation}
The gradient of \eqref{eq:local:energy} with respect to $x_{=k}$ is zero when
\begin{equation}
	\label{eq:local:system}
	(\tens{x}_{\neq k}^* \tens{K} \tens{x}_{\neq k}) x_{=k} 
	= 
	\tens{x}_{\neq k}^* \tens{b}
	.
\end{equation}
The problem size of this reduced problem is small and we can use standard numerical methods such as a direct solver to solve this subproblem efficiently. 
For the general ALS approach the resulting TT-ranks --- and therefore the maximum accuracy --- are fixed by the ranks set in the initial guess.
But with some extension we can adapt the TT-ranks of the solution dynamically. 
For the purpose of calculating a low-rank approximation, we chose the so-called \textit{Alternating Minimal Energy} (AMEn) method \cite{DolgovSavostyanov:2014:1}.

The AMEn method proved to be robust and has a fast convergence rate. 
For a detailed analysis and more information we refer the reader to \cite{DolgovSavostyanov:2014:1}. 
Note that we can also use this method to compute a matrix-matrix multiplication and it is especially efficient when using a large number of time steps and/or large number of output points as the memory reduction gets more significant in these cases.

Based on this discussion, matrix-vector products with the prior-pre\-con\-ditio\-ned misfit Hessian $\gPrior \hMisfit$ or its equivalent analog \eqref{eq:hMis:genEigO} 
\begin{equation*}
	\begin{aligned}
		% \gPrior \hMisfit 
		% &
		% = 
		% \gPrior(\tens{M}_0 \tens{K}^{-\transp} \tens{C}^\transp) \gNoise^{-1} (\tens{C} \tens{K}^{-1} \tens{M}_0) 
		% \\
		\gPriorO^{T/2} \hMisfit \gPriorO^{1/2}
		&
		=
		\gPriorO^{T/2} (\tens{M}_0 \tens{K}^{-\transp} \tens{C}^\transp) \gNoise^{-1} (\tens{C} \tens{K}^{-1} \tens{M}_0) \gPriorO^{1/2}
	\end{aligned}
\end{equation*}
can be fully performed within the tensor-train format.
Notice that the action of $\gPriorO^{1/2}$ and its transpose are available as $\gPriorO^{1/2} = A_h^{-1} M$.
The numerical implementation uses this alternative formulation of the generalized eigenvalue problem \eqref{eq:hMis:genEigO}.

\section{Numerical Results}
\label{sec:numerical_results}

This section is devoted to the presentation of numerical results and demonstration of efficiency of the tensor train approximation of the posterior covariance.
All experiments are based on the machine described in \cref{sec:the_thermo_elastic_problem} with the physical and numerical parameters as in \cref{table:setupPO}. 
We provide the data and the code in a container image \cite{NaumannBuenger:2023:1}.

In both approaches the inverse of $A_h$ in the prior was approximated with a (single) sparse LU decomposition, based on \umfpack \cite{Davis:2004:2}.
For the computation of the sensitivity matrix~$F$ \eqref{eq:derivative:continous} in the direct approach, we also use a single sparse LU decomposition using \umfpack of the constant system matrix in \eqref{eq:fully_discrete_adjoint_problem}.
In the direct approach we used \scipy (\cite{VirtanenGommersOliphantHaberlandReddyCournapeauBurovskiPetersonWeckesserBrightVanDerWaltBrettWilsonMillmanMayorovNelsonJonesKernLarsonCareyPolatFengMooreVanderPlasLaxaldePerktoldCimrmanHenriksenQuinteroHarrisArchibaldRibeiroPedregosaVanMulbregt:2020:1}) for the solution of the generalized eigenvalue problem \eqref{eq:hMis:genEig} and the subsequent evaluation of the approximate posterior variance \eqref{eq:evaluation_of_diagonal_entries_of_posterior_covariance}.
The tensor train approach was implemented in \matlab, using the TT~toolbox (\url{https://github.com/oseledets/TT-Toolbox}).

\begin{table}
	\begin{footnotesize}
		\centering
		\begin{tabular}{l | r | l}
			parameter & value & description \\
			\hline
			$n_T$ & $120$ & number of time steps \\
			$T$ & $\SI{120}{\second}$ & final time \\
			$a^{(i)}$ &(\num[round-precision=4]{0.0618533085154064},\num[round-precision=4]{0.2499971726638156},\num[round-precision=4]{0.13053275665353803},\num[round-precision=4]{0.25743961872851184},\num[round-precision=4]{0.1402166704168589}) & parameters for prior variance \\
			$b^{(i)}$ & (\num[round-precision=4]{1.9350600255384578},\num[round-precision=4]{7.821077755264953},\num[round-precision=4]{4.083673541257472},\num[round-precision=4]{8.053912185915216},\num[round-precision=4]{4.386631537586745}) & parameters for prior variance \\
			$\tau$ & \num{1800} & parameters for prior variance \\
			$\sigmanoise$ & \SI[round-precision=1]{0.1}{\kelvin} & measurement standard deviation \\
			$r$ & $50$ & rank of the posterior approximation \\ 
			% tol & $10^{-5}$ & tolerance for the AMEn solver \\
			% truncation tol & $10^{-7}$ & truncation tolerance
		\end{tabular}
		\caption{Setup for the numerical experiments. For the construction of the prior variance, see \cref{sec:prior_selection}. For the posterior approximation, see \cref{subsection:Low-rank_approximation}.}
		\label{table:setupPO}
	\end{footnotesize}
\end{table}

\begin{figure}
  \begin{minipage}[t]{0.48\textwidth}
    \includegraphics[width=\textwidth]{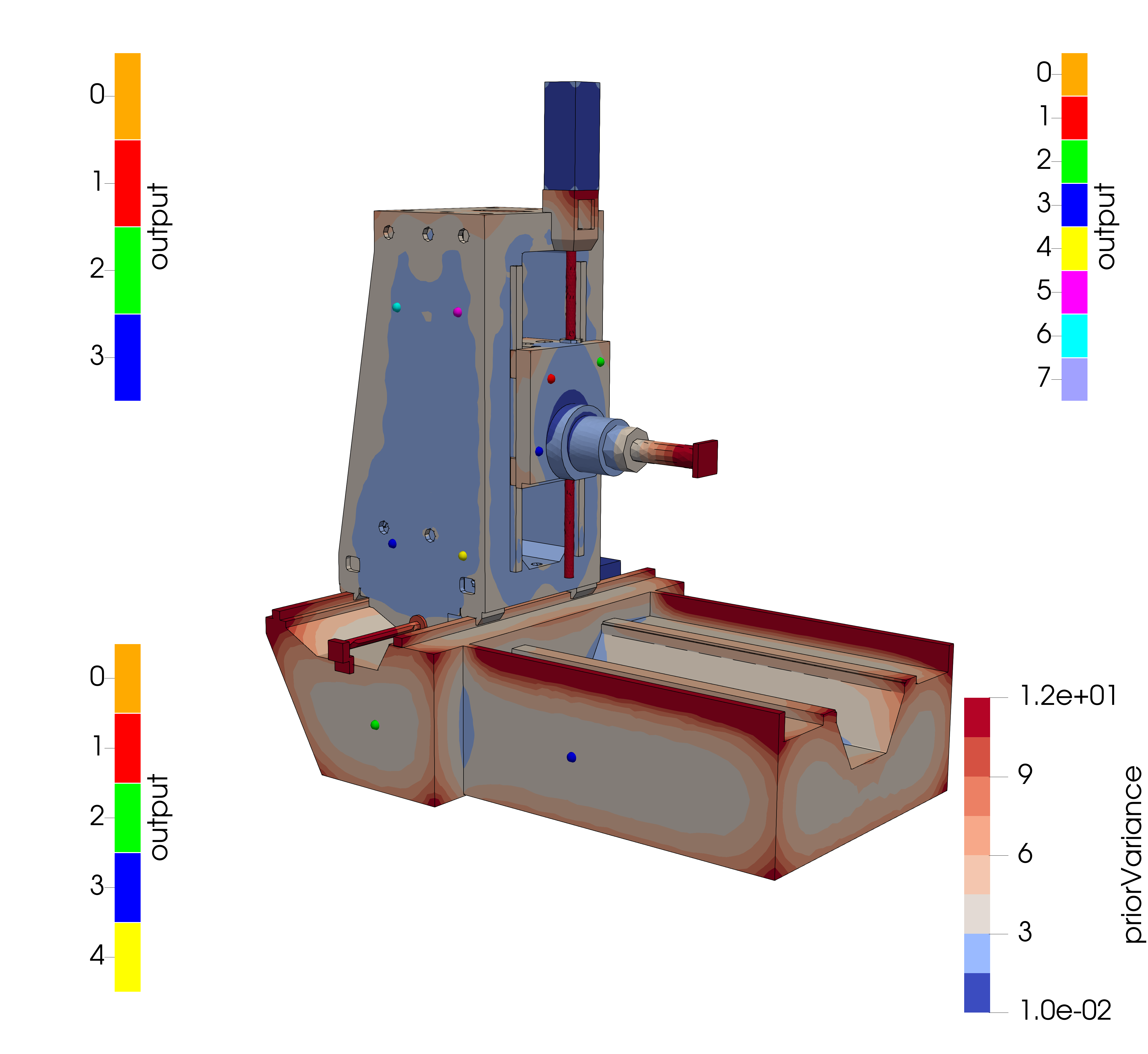}
    \caption{Prior variance field using the prior operator \eqref{eq:prior:discretePDE} with the parameters from \cref{table:setupPO}.}%
    \label{fig:priorVariance:physicalParameters}
  \end{minipage}
	\hfill
  \begin{minipage}[t]{0.48\textwidth}
    \includegraphics[width=\textwidth]{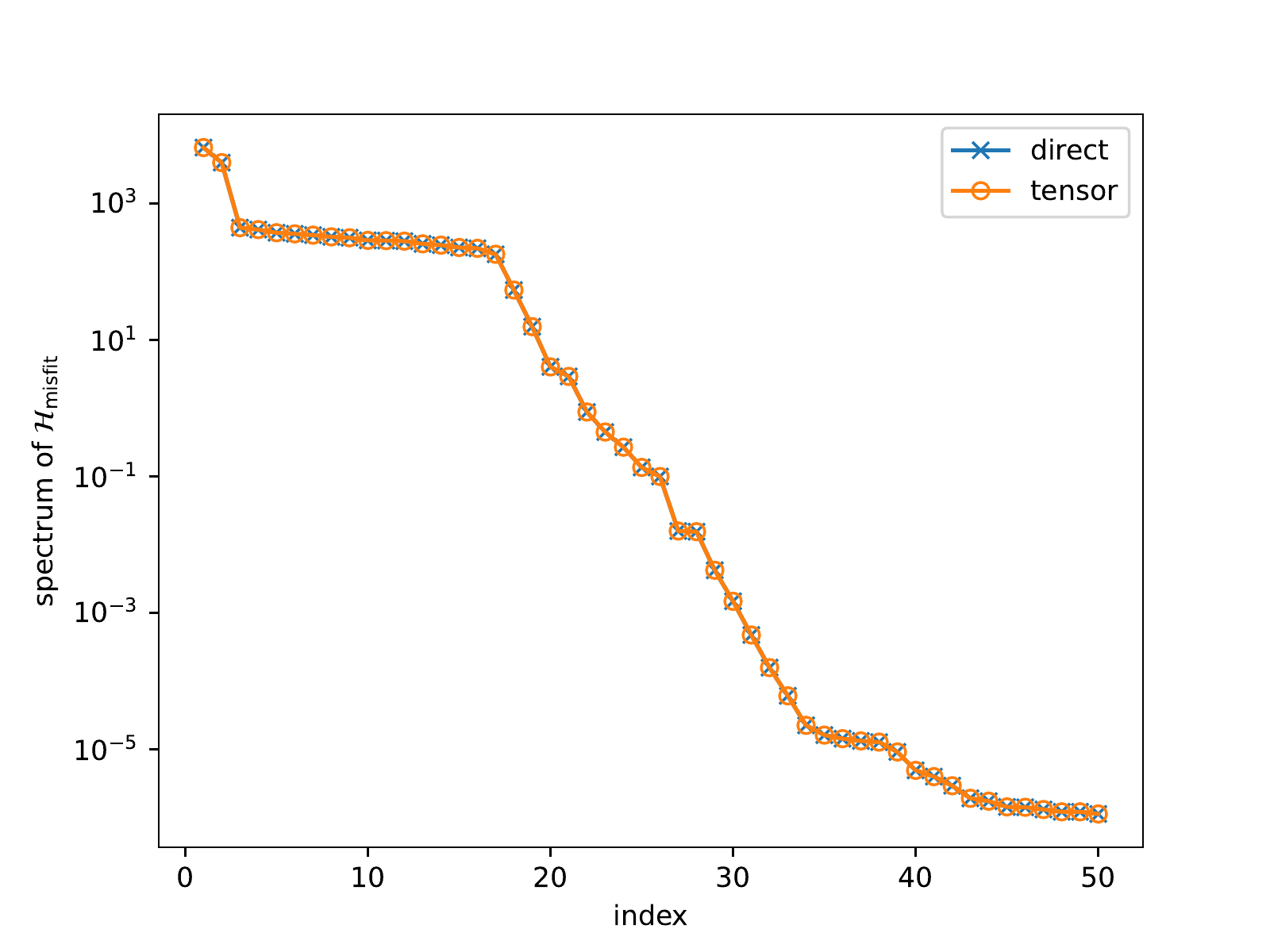}
		\caption{Spectrum \eqref{eq:hMis:genEig} of the generalized eigenvalue problem $(\hMisfit, \gPrior^{-1})$ using the direct and the tensor train approach.} 
    \label{figure:SVD:compare}
  \end{minipage}
\end{figure}

We begin with the presentation of the prior variance field in \cref{fig:priorVariance:physicalParameters}.
The color scale ranges from \SI[round-precision=2]{0.01}{\kelvin\squared} to \SI{12}{\kelvin\squared} in steps of \SI[round-precision=2]{1.5}{\kelvin\squared}.
This color scale is chosen based on the values of the \emph{posterior} variance discussed below.
The largest values of the prior variance occur near the machine's edges.
These maximal values reach up to \SI{35}{\kelvin\squared}, which corresponds to a standard deviation of around \SI[round-precision=1]{5.9}{\kelvin}.
On the surface, far away from the edges, the prior variance is between the average prior variance \SI{3}{\kelvin\squared} and \SI{9}{\kelvin\squared}.
The minimal value of the prior variance is about \SI[round-precision=1]{1.7}{\kelvin\squared}.
Since the point values of the prior variance are sufficiently larger than the measurement variance $\sigmanoise^2$, we can expect the information gained by the measurements to be significant and the posterior variance to be notably different from the prior variance.
Notice that \cref{fig:priorVariance:physicalParameters} also shows the location of the $17$~temperature sensors which were already shown in \cref{fig:ACW:Outputs}.

In \cref{figure:SVD:compare} we depict the leading $r = 50$ spectral values of the eigenvalue problems \eqref{eq:HmisApprox} and \eqref{eq:hMis:genEigO} related to the data misfit Hessian, once using the direct approximation approach of \cref{sub:direct_approximation}, and once using the tensor train approximation of \cref{sub:low_rank_time_discretization}.
Both approaches are in excellent agreement over a span of 8~orders of magnitude, down to a value of at least \num{1.0e-5}.

\begin{figure}
  \begin{minipage}[t]{0.48\textwidth}
    \includegraphics[width=0.95\textwidth]{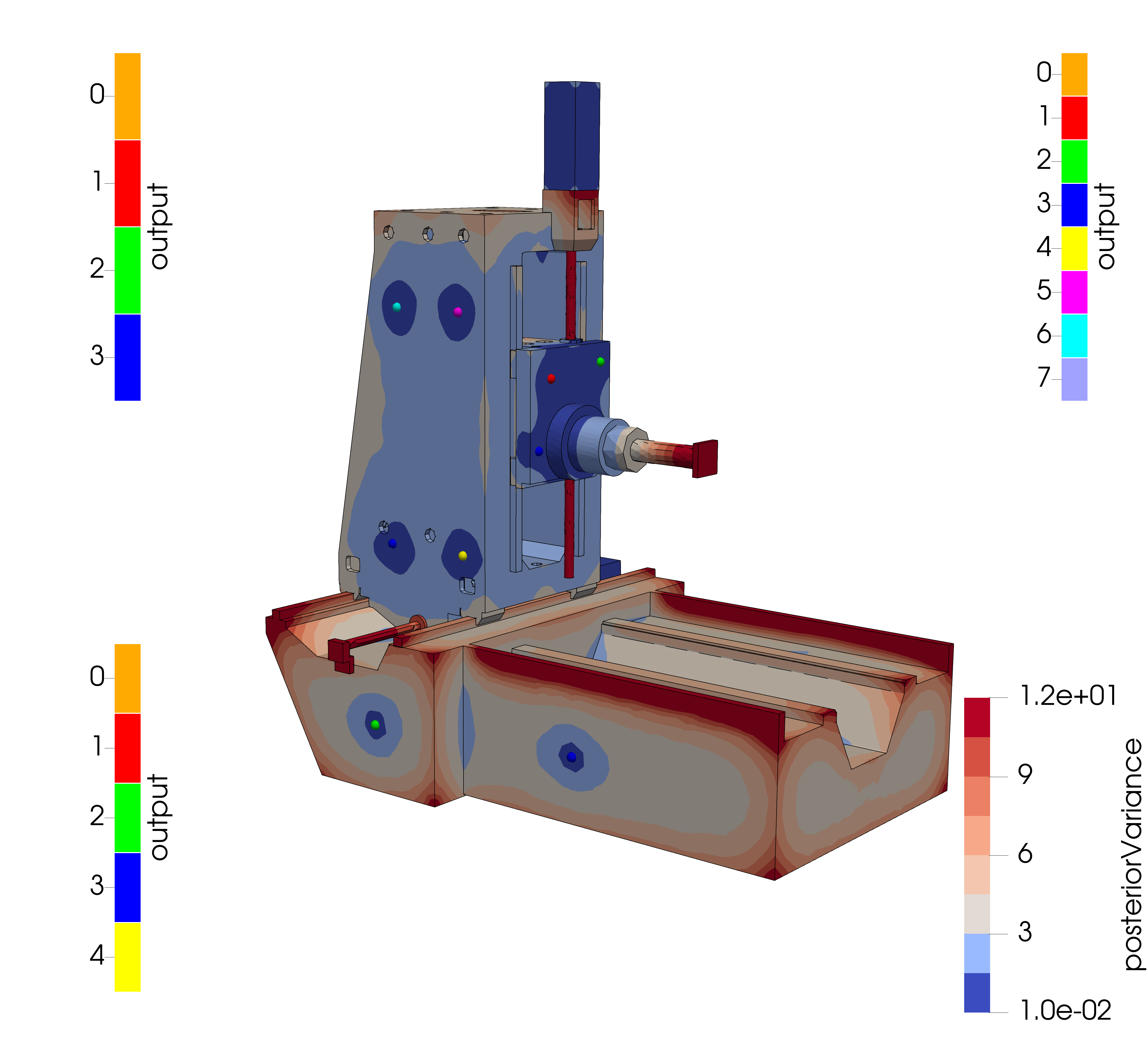}
    \caption{Posterior variance \eqref{eq:evaluation_of_diagonal_entries_of_posterior_covarianceO}, approximated with $39$~eigenpairs using the tensor train approach.}
    \label{figure:posteriorVariance}
  \end{minipage}
	\hfill
  \begin{minipage}[t]{0.48\textwidth}
     \includegraphics[width=0.95\linewidth]{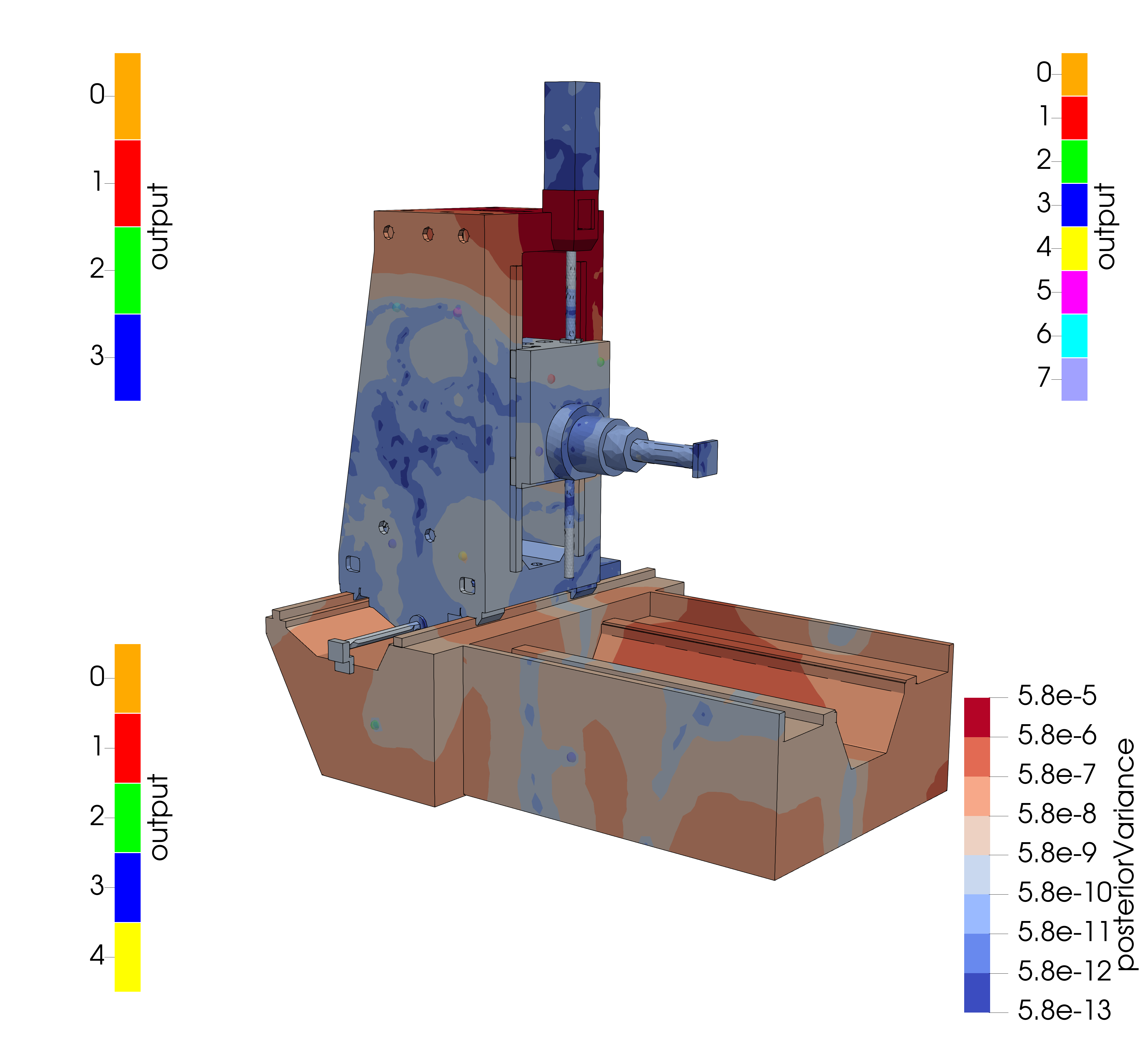}
		\caption{Absolute difference in the posterior variance between using the tensor train approach \eqref{eq:evaluation_of_diagonal_entries_of_posterior_covarianceO} and the direct approach \eqref{eq:evaluation_of_diagonal_entries_of_posterior_covariance}.}%
    \label{fig:diff:tensorDirect}
  \end{minipage}
\end{figure}

\Cref{figure:posteriorVariance} depicts the values of the posterior variance.
Its nodal values are evaluated using \eqref{eq:evaluation_of_diagonal_entries_of_posterior_covarianceO}, using the leading $r = 39$ eigenpairs and the tensor train approximation of the misfit Hessian.
We see that the covariance is lowest next to the measurements in the range of the measurement variance, as expected.
Although \cref{figure:SVD:compare} confirms that the generalized eigenvalues of the misfit Hessian agree even beyond index~$39$, this is not necessarily the case also for the eigenvectors.
We therefore display in \cref{fig:diff:tensorDirect} the absolute difference of the posterior variance using both the tensor train and the direct approach.
The latter also used $39$~eigenpairs.
The absolute difference is below \SI{6e-5}{\kelvin\squared} and thus well below the order of the measurement variance $\sigmanoise^2$.
This confirms that from a practical point of view, the tensor train approximation of the misfit Hessian does not contribute significantly to the posterior uncertainty in the initial temperature.

\begin{table}
	\begin{footnotesize}
		\centering
		\begin{tabular}{l | r | r | r |}
			& \python & \matlab & \matlab \\
			& direct ($r = 50$) & tensor ($r = 26$) & tensor ($r = 39$) \\ \hline
			prior variance (offline phase) & \num{682} & \num{1991} & \num{1991} \\ \hline
			sensitivity matrix & \num{238} &	& \\
			generalized eigenvalue problem & \num{193} & \num{3271} & \num{6305}\\
			posterior variance & \num{58}	& \num{56}	& \num{94}	\\ \hline
			total & \num{489} & \num{3327}	& \num{6399}
		\end{tabular}
		\caption{Runtimes in seconds for the computation of the sensitivity matrix \eqref{eq:derivative:continous}, solution of the generalized eigenvalue problem \eqref{eq:hMis:genEig} (direct approach) and \eqref{eq:hMis:genEigO} (tensor train approach, for $r = 26$ and $r = 39$ leading values), and approximation of the posterior variance \eqref{eq:evaluation_of_diagonal_entries_of_posterior_covariance} (direct approach) and \eqref{eq:evaluation_of_diagonal_entries_of_posterior_covarianceO} (tensor train approach).}
		\label{table:runtimes}
	\end{footnotesize}
\end{table}

%\added[id=RH]{\textbf{@AN: Comment that the prior variance does not need to be recomputed, \eg, when the sensor location changes.}}

The experiments were run on an Intel i7 CPU with \SI{8}{\mebi\byte} cache and \SI{16}{\gibi\byte} RAM using four cores.
In \cref{table:runtimes} we denote the runtime of each sub task for both approaches.
The runtimes for the evaluation of the prior variance $e_k^\transp \gPrior \, e_k$, \ie, the part of the posterior variance attributed to the prior, differ due to different implementations in \python and \matlab, respectively, which also use different underlying \blas implementations.
While the differences in runtime are significant, the evaluation of the prior variance can be considered an offline part of any computation since it only depends on the geometry and, particular, does not depend on sensor positions. 

\begin{figure}
	\centering
	\includegraphics[width=0.5\textwidth]{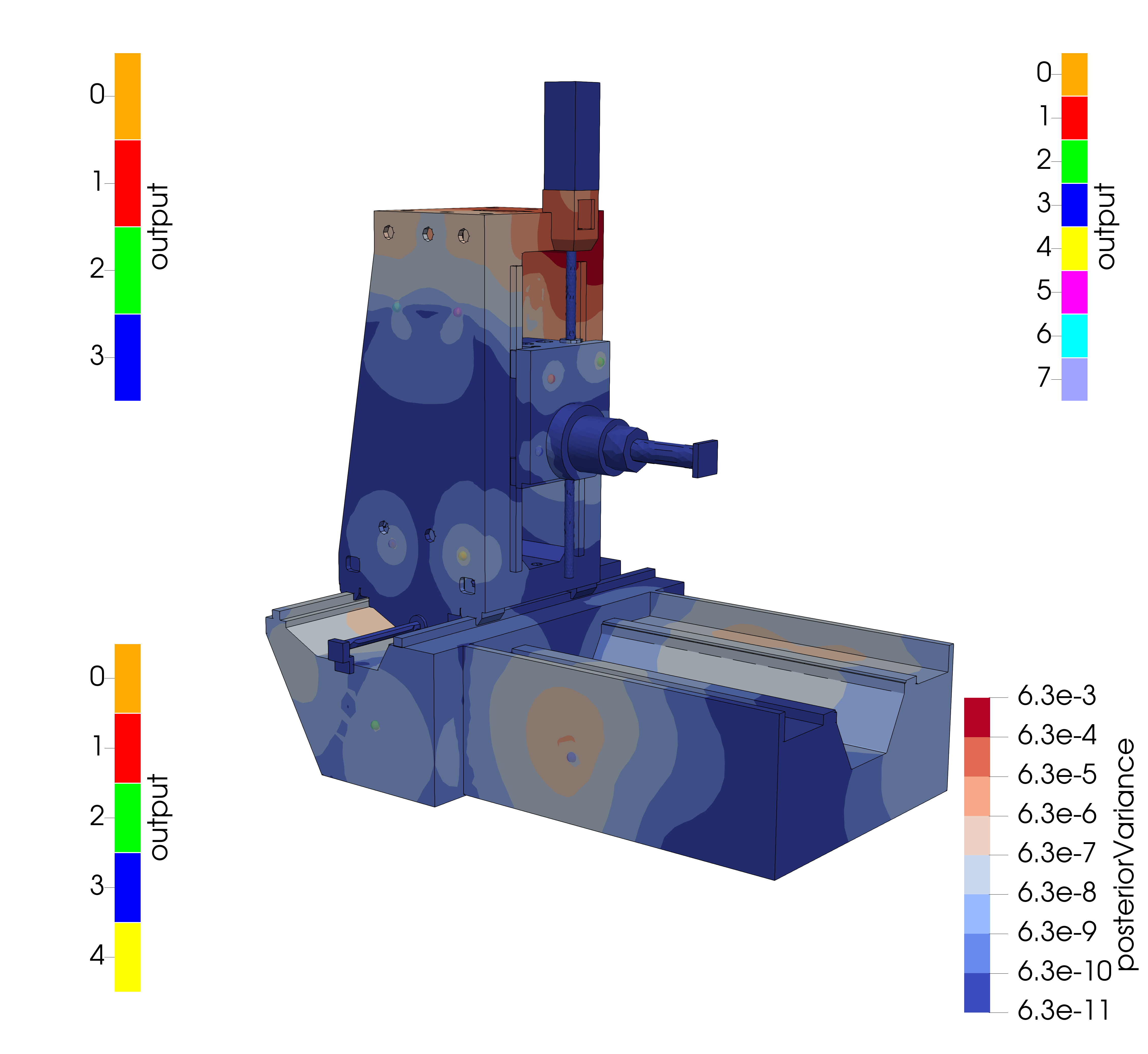}
	\caption{Relative difference of the posterior variance using the tensor train approach between $r = 26$~eigenvalues relative to using $r = 39$ eigenvalues.}
	\label{figure:diffLowRank}
\end{figure}

The direct approach spent about one half of the time for the computation of the sensitivity matrix and the other half for the solution of the generalized eigenvalue problem.
The final evaluation of the posterior variance \eqref{eq:evaluation_of_diagonal_entries_of_posterior_covariance} is not significant in comparison.
By contrast, the evaluation of the posterior variance \eqref{eq:evaluation_of_diagonal_entries_of_posterior_covarianceO} using the tensor train approach is dominated by the runtime to solve the generalized eigenvalue problem.
This can be attributed to the fact that the repeated \enquote{matrix-vector} multiplications with $\hMisfit$ require repeated evaluations of \eqref{eq:parameter_to_observables}.
Although the latter uses tensor train approximations, it is still costly compared to the once-and-for-all evaluation of the sensitivity matrix.
Reducing the number of eigenvalues from $r = 39$ to $r = 26$, which offers sufficient accuracy as shown in \cref{figure:diffLowRank}, cuts the runtime down to about half.
The observed runtime thus seems to increase slightly stronger than linear with~$r$, which we attribute to the restart of the Arnoldi method.

The major disadvantage of the direct approach is the memory consumption of the sensitivity matrix~$F$ \eqref{eq:derivative:continous}.
Its dimension is $n_t \, n_y$ (number of observation times multiplied with the number of outputs or sensors) times $n_x$ (number of degrees of freedom in the discretization).
For the problem at hand, we recall $n_t = 120$, $n_y = 17$ and $n_x$ is approximately~$\num{75000}$.
Even for this moderate number of outputs, the sensitivity matrix thus consumes approximately \SI[round-precision=1]{1.2}{\gibi\byte} of memory.
For applications such as optimized sensor placement, where the number of outputs $n_y$ is in the order of $n_x$, the memory consumption of the direct approach is prohibitive.
By contrast, the dominant storage in the tensor train approach are the coefficient matrices $\phys{K}$ and $\phys{M}$, which are needed in any case.

\section{Conclusion and Outlook}
\label{sec:conc}

We derived a coupled thermal finite element model for the simulation of the temperature fields in a machine tool.
One characteristic of this model is its large dimension of the state space versus a small number of outputs, in our case temperature sensors.
We considered the inverse problem of estimating the initial temperature field using time series of temperature measurements.

Using a Bayesian approach, our main interest lies with the evaluation of the posterior initial temperature variance (as a function of space).
To make this computationally tractable, we approximate the data misfit Hessian by a low-rank approximation using its leading eigenpairs.
The iterative solution of this eigenvalue problem requires the repeated evaluation of matrix-vector products with the linearized parameter-to-observable map and its transpose.
We compare two approaches to numerically realize these matrix-vector products.
In the direct approach, we set up the linearized parameter-to-observable map as a full matrix.
The major disadvantage of the direct approach is its memory consumption.
Therefore, we consider a tensor train approach as an alternative, which does not need significantly more memory than that for the finite element matrices.

While the tensor train approach turned out to require more CPU time in our implementation, the memory footprint of the direct approach quickly becomes intractable as the number of outputs increases. 
Here, low-rank tensor formats typically alleviate these effects as they scale linearly in the number of outputs. 
We will investigate this further especially when optimal sensor placement tasks are considered.	
Moreover, we can anticipate that for \emph{nonlinear} parameter-to-observable maps, the direct approach becomes prohihitive for even smaller problem sizes and nonlinear space-time approach in low-rank tensor form could be applicable here \cite{DolgovStoll:2017:1}.

% Insert the appendix.
\appendix

% Insert the bibliography.
\printbibliography

\end{document}